\documentclass{amsart}

\usepackage{amssymb}
\usepackage[hidelinks]{hyperref}
\usepackage{verbatim} % For commenting.
\usepackage{xypic} % For commutative diagrams.
\usepackage[mathcal]{euscript} % Change font for mathcal.
\usepackage{array} % for \newcolumntype macro
\numberwithin{equation}{subsection}
\newtheorem{prop}[equation]{Proposition}
\newtheorem{thm}[equation]{Theorem}
\newtheorem{lem}[equation]{Lemma}
\newtheorem{cor}[equation]{Corollary}
\theoremstyle{remark}
\newtheorem{remark}[equation]{Remark}
\theoremstyle{definition}
\newtheorem{defn}[equation]{Definition}
\newtheorem{example}[equation]{Example}

\newcolumntype{C}{>{$}c<{$}}
\renewcommand{\lto}{\leftarrow}
\newcommand{\bbz}{\mathbb{Z}}
\newcommand{\bbn}{\mathbb{N}}
\newcommand{\mfm}{\mathfrak{m}}
\newcommand{\dev}{\varepsilon}
\DeclareMathOperator{\upB}{\tn{B}}
\DeclareMathOperator{\im}{im}

\DeclareMathOperator{\upP}{\tn{P}}
\DeclareMathOperator{\upZ}{\tn{Z}}
\DeclareMathOperator{\upH}{H}
\newcommand{\sus}{{\scriptstyle\Sigma}}
\DeclareMathOperator{\bideg}{bideg}
\DeclareMathOperator{\cls}{cls}
\newcommand{\modu}[1]{\ \tn{(mod $#1$)}}
\DeclareMathOperator{\supp}{supp}

\DeclareMathOperator{\depth}{depth}
\DeclareMathOperator{\sfm}{\mathsf{M}}
\newcommand{\dotimes}{\otimes^\mathsf{L}}
\newcommand{\op}{\mathsf{op}}
\newcommand{\surj}{\twoheadrightarrow}
\newcommand{\tn}[1]{\textnormal{#1}}
\newcommand{\cl}[1]{\overline{#1}}
\newcommand{\ang}[1]{\left\langle #1 \right\rangle}
\newcommand{\bd}{\partial}
\DeclareMathOperator{\Tor}{Tor}

\DeclareMathOperator{\Hom}{Hom}

\newcommand{\mbe}{\mathbf{e}}
\newcommand{\mbp}{\mathbf{p}}
\newcommand{\mbu}{\mathbf{u}}
\newcommand{\mbv}{\mathbf{v}}
\newcommand{\mbw}{\mathbf{w}}

\title[Linear resolutions over Koszul complexes]{Linear resolutions over Koszul complexes and Koszul homology algebras}
\author{John Myers}
\address{Department of Mathematics, SUNY Oswego, Oswego, New York, USA}
\email{john.myers@oswego.edu}


\begin{document}

\begin{abstract}
Let $R$ be a standard graded commutative algebra over a field $k$, let $K$ be its Koszul complex viewed as a differential graded $k$-algebra, and let $H$ be the homology algebra of $K$. This paper studies the interplay between homological properties of the three algebras $R$, $K$, and $H$. In particular, we introduce two definitions of Koszulness that extend the familiar property originally introduced by Priddy: one which applies to $K$ (and, more generally, to any connected differential graded $k$-algebra) and the other, called \textit{strand-Koszulness}, which applies to $H$. The main theoretical result is a complete description of how these Koszul properties of $R$, $K$, and $H$ are related to each other. This result shows that strand-Koszulness of $H$ is stronger than Koszulness of $R$, and we include examples of classes of algebras which have Koszul homology algebras that are strand-Koszul.
\end{abstract}

\maketitle
\thispagestyle{empty}

\section*{Introduction}

%\begin{comment}
Koszul complexes are classical objects of study in commutative algebra. Their structure reflects many important properties of the rings over which they are defined, and these reflections are often encoded in the product structure of their homology. Indeed, Koszul complexes are more than just complexes --- they are, in fact, the prototypical examples of differential graded ($=$ DG) algebras in commutative ring theory, and thus the homology of a Koszul complex carries an algebra structure. These homology algebras encode (among other things) the Gorenstein condition \cite{AvramovGolod1971}, the Golod condition \cite{Golod1962}, and whether or not the ring is a complete intersection \cite{Assmus1959}. The present paper follows in the spirit of these results by studying how certain homological properties of Koszul complexes and their homology algebras are interrelated with homological properties of the rings from which they originate.

In more detail, let $R$ be a standard graded commutative algebra over a field $k$, let $K$ be a Koszul complex on a minimal set of generators of the augmentation ideal of $R$ (which we call \textit{the} Koszul complex of $R$), and let $H$ be the homology algebra of $K$. The ground field $k$ can then be resolved by free modules over the algebras $R$ and $H$, and by semifree modules over the algebra $K$ (the definition of \textit{semifree module} belongs to \textit{DG homological algebra}, and it will be recalled in the body of the paper). Our goal is then to compare the three resolutions of $k$ over the three different algebras and to identify when linearity of one of the resolutions implies linearity of one of the others. Linearity of the resolution over $R$ means that the algebra $R$ belongs to the class of \textit{Koszul algebras}, a familiar class whose theory was first systematically laid down by Priddy in \cite{Priddy1970} and which has become a class of central importance in the homology theory of both noncommutative and commutative graded algebras. See, for example, the survey \cite{Froberg1999} and the monograph \cite{PolishchukPositselski2005} for expositions on the noncommutative side, and the survey \cite{Conca2013} directed at applications in commutative algebra.

On the other hand, what we mean by ``linearity'' of the resolutions over $K$ and $H$ requires further explanation. As we will recall below, there is a DG version of $\Tor^K{(k,k)}$ which can be computed via a semifree resolution of $k$ over $K$. The grading on $R$ induces an extra internal degree on $\Tor^K{(k,k)}$ in addition to its usual homological degree, and we then say that $K$ is a \textit{Koszul DG algebra} when the ``non-diagonal'' part of $\Tor^K{(k,k)}$ vanishes, in complete analogy with the classical definition of Koszulness.

This definition follows in the spirit of (but is not exactly identical to) the definition offered by He and Wu in \cite{HeWu2008} who deal with connected \textit{cochain} DG algebras instead of \textit{chain} algebras, which is where our focus lies. The differential accompanying the former type of algebra raises homological degree, while the one accompanying the latter lowers homological degree. %Cochain DG algebras tend to occur more naturally in geometric and topological settings whereas chain DG algebras occur in algebraic ones; for example, He and Wu's first example of a Koszul cochain DG algebra is the minimal model of the de Rham complex of a smooth manifold satisfying certain extra hypotheses. See also the dictionary \cite{AvramovHalperin1986}, where some distinctions (and similarities) between the topological and algebraic sides of DG homological algebra are elucidated.

The relationship between Koszulness of the algebras $R$ and $K$ is then satisfyingly simple and should not surprise specialists in this field.

\theoremstyle{theorem}
\newtheorem*{thmA}{Theorem A}
\begin{thmA}
Let $R$ be a standard graded commutative algebra over a field $k$ and let $K$ be its Koszul complex. The algebras $R$ and $K$ are Koszul simultaneously.
\end{thmA}

%With all of the definitions securely in place, the proof comes by repackaging a construction due to Tate. In \cite{Avramov1984}, Avramov wrote that ``...for most purposes one can replace the ring by its Koszul complex''; our Theorem A lends more evidence to the veracity of this claim.

Shifting attention to the homology algebra $H$, the grading on $R$ also induces an internal degree on $H$ giving it the structure of bigraded $k$-algebra. It is easy to show that $H_{ij}\neq 0$ implies $i=j=0$ or $j-i>0$, a property that we describe concisely by saying $H$ \textit{lives in positive strands}. This shape makes it possible to perform a totalization process, which we call \textit{strand totalization}, to produce a connected graded $k$-algebra $H'$ with $H'_n = \bigoplus_{j-i=n} H_{ij}$ for each $n$. Then, when we speak of $k$ having a linear free resolution over $H$, we mean that $H'$ is a Koszul algebra in the classical sense, in which case we say $H$ is \textit{strand-Koszul}.

The complete relationship between the three notions of Koszulness is expressed in

\theoremstyle{theorem}
\newtheorem*{thmB}{Theorem B}
\begin{thmB}
Let $R$ be a standard graded commutative algebra over a field $k$, let $K$ be its Koszul complex, and let $H$ be the homology algebra of $K$. The following are equivalent.
\begin{enumerate}
\item The algebra $H$ is strand-Koszul.
\item The algebra $K$ is Koszul and quasi-formal.
\item The algebra $R$ is Koszul and $K$ is quasi-formal.
\item The algebra $R$ is Koszul and the minimal graded free resolution of $k$ over $R$ has an Eilenberg-Moore filtration.
\end{enumerate}
\end{thmB}

\textit{Quasi-formality} and the existence of an \textit{Eilenberg-Moore filtration} are two technical conditions which facilitate a clean transfer of homological properties from $R$ and $K$ to $H$. The former condition, introduced and studied by Positselski in \cite{Positselski2017}, concerns the degeneracy of a spectral sequence which connects $H$ to $K$, while an Eilenberg-Moore filtration allows one to ``build'' the minimal free resolution of $k$ over $H$ starting from the minimal free resolution of $k$ over $R$.

Theorem B thus shows that strand-Koszulness of $H$ is stronger than Koszulness of $R$ and $K$. For example, the Koszul algebra
	\[R = \frac{k[x,y,z,u]}{(x^2,xy,xz+u^2,xu,y^2+z^2,zu)},
	\]
which is listed as ``isotope 63ne'' in Roos' catalog \cite[Appendix A]{Roos2016}, has a Koszul homology algebra which is not strand-Koszul (verification by \texttt{Macaulay2} \cite{Macaulay2} using the \texttt{DGAlgebras} package written by Frank Moore). We note that out of the 104 algebras of embedding dimension $4$ that Roos catalogs, there are only two Koszul algebras which have Koszul homology algebras that are not strand-Koszul. Thus, in combination with our next result, we have that strand-Koszulness of $H$ is ``nearly'' equivalent to Koszulness of $R$, as long as $R$ has embedding dimension $\leq 4$.

\theoremstyle{theorem}
\newtheorem*{thmC}{Theorem C}
\begin{thmC}
Let $R$ be a standard graded commutative algebra over a field $k$, let $H$ be the Koszul homology algebra of $R$, and suppose that one of the following statements is true.
\begin{enumerate}
\item The algebra $R$ is Koszul with embedding dimension $\leq 3$.
\item The algebra $R$ is Koszul and Golod.
\item The algebra $R$ is a quadratic complete intersection.
\item The algebra $R$ is artinian Gorenstein of socle degree $2$, and either $k$ does not have characteristic $2$ or the embedding dimension of $R$ is odd.
\item The algebra $R$ is Koszul and the defining ideal of $R$ is minimally generated by three elements.
\item The defining ideal of $R$ is the edge ideal of a path on $\geq 3$ vertices.
\end{enumerate}
Then $H$ is strand-Koszul.
\end{thmC}

The proof uses a mixture of techniques. Previously established theory makes quick work of the implications ``(1)-(3) $\Rightarrow$ $H$ is strand-Koszul,'' while our proof of the implications ``(4)-(6) $\Rightarrow$ $H$ is strand-Koszul'' is based on (noncommutative) Gr\"obner basis techniques and computations of presentations of $H$ as a quotient of a free algebra.

We point out that the strand totalizations $H'$ are not new, except perhaps for the name. For example, a question posed by Boocher, D'Al\`\i, Grifo, Monta\~no, and Sammartano \cite{BoocherEtAl2015} asks whether or not the algebra $H'$ is generated by degree-$1$ elements when $R$ is a Koszul integral domain (this replaces the same question posed earlier by Avramov without the hypothesis that $R$ is an integral domain --- the answer to that question turned out to be negative, see \cite[Remark 3.2]{BoocherEtAl2015}). In \cite{FrobergLofwall2002}, Fr\"oberg and L\"ofwal study $H'$ and uncover a connection between this algebra and (part of) the homotopy Lie algebra of $R$. The relation between strand totalizations and Koszulness of $R$ was also considered by Croll et al.\ \cite{CrollEtAl2016} who showed that if there is an element $\zeta\in H_{1,2} \subseteq H'_1$ such that $H'\cdot \zeta = H'_{>1}$, then $R$ is Koszul.

Section 1 of this paper sets up definitions and notations and proves some preliminary results, while section 2 is dedicated to the proofs of Theorems A and B and also contains some results on relationships between Poincar\'e-Betti series, Hilbert series, and low-degree Betti numbers. The third and final section contains the proof of Theorem C.
%\end{comment}

\section{Definitions and preliminaries}

Throughout this paper $k$ denotes a field. Elements of graded objects will always be assumed homogeneous and thus all elements of such objects have degrees.

We assume that the reader is familiar with the basic terminology and theory of connected $\bbz$-graded $k$-algebras; everything that we need can be found in sections 0 and 4 in Chapter 1 of \cite{PolishchukPositselski2005}, or the Appendice in \cite{Lemaire1974}. We single out just a few definitions and terms from this theory that are of particular importance to us. First, we make

\begin{defn}\label{defn:Koszul}
Let $A$ be a connected $\bbz$-graded $k$-algebra. We will say $A$ is \textit{Koszul} if $\beta_{ij}^A(k) \neq 0$ implies $i=j$ for all $i,j\geq 0$.
\end{defn}

Here
	\[\beta_{ij}^A(k) = \dim_k{ \Tor^A_{i}{(k,k)}_j}
	\]
is the \textit{$(i,j)$-th Betti number} of $k$. All of these numbers are finite, and the generating function built from them is denoted
	\[\upP^A_k(s,t) = \sum_{i,j} \beta^A_{ij}(k) s^i t^j
	\]
and is called the \textit{bigraded Poincar\'e series} of $k$.

A $\bbz$-graded module $M$ over a connected $\bbn$-graded $k$-algebra is said to be \textit{finite} if $M_j=0$ for all $j\ll 0$ and each $M_j$ is finite-dimensional over $k$. The \textit{$q$-th shift} of $M$ is defined to be the $\bbz$-graded module
	\[M(q) = \bigoplus_{j\in \bbz} M(q)_j \quad \tn{with} \quad M(q)_j = M_{q+j}.
	\]

\subsection{Connected $\bbz^2$-graded algebras}\label{subsec:bigraded}

This paper is concerned with both $\bbz$- and $\bbz^2$-graded objects; hereafter, the former objects will simply be called \textit{graded}, while the latter objects will be called \textit{bigraded}.

A bigraded $k$-algebra $A = \bigoplus_{i,j\in \bbz} A_{ij}$ is said to be \textit{connected} when $A_{00} = k$; $A_{ij}\neq 0$ implies $i,j\geq 0$; and each homogeneous component $A_{ij}$ is finite-dimensional. Given an element $a\in A_{ij}$, the integers $i$ and $j$, the pair $(i,j)$, and the sum $i+j$ will be called, respectively, the \textit{homological degree}, \textit{internal degree}, \textit{bidegree}, and \textit{total degree} of $a$. The homological degree of $a$ will be denoted $|a|$, its internal degree denoted $\deg{a}$, and its bidegree denoted $\bideg{a}$. We will write $A_i$ for the graded $k$-space $A_{i,\ast}$.

Let $A$ be a connected bigraded $k$-algebra and $M$ a bigraded left $A$-module.

The \textit{$p$-th homological shift} of $M$ is the bigraded module
	\[\sus^p M = \bigoplus_{i,j\in \bbz} (\sus^p M)_{ij}, \quad \tn{with} \quad (\sus^p M)_{ij} = M_{i-p,j}.
	\]
The \textit{$q$-th internal shift} of $M$ is the bigraded module
	\[M(q) = \bigoplus_{i,j\in \bbz} M(q)_{ij}, \quad \tn{with} \quad M(q)_{ij} = M_{i,q+j}.
	\]

The module $M$ is said to be \textit{finite} if $M_i=0$ for all $i\ll 0$ and each $M_i$ is finite as an $A_0$-module (as above, we write $M_i$ for $M_{i,\ast}$).

The field $k$ can be viewed as a bigraded $k$-algebra concentrated in bidegree $(0,0)$. The natural map $A \to k$ (called the \textit{augmentation}) is a morphism of bigraded algebras and gives $k$ a bigraded $A$-module structure. We define the \textit{augmentation ideal} of $A$, denoted $\mfm_A$, to be the ideal $\ker{(A \to k)}$. If the algebra $A$ is understood, we will write $\mfm$ in place of $\mfm_A$.

If $M$ is finite, then $M$ has a bigraded minimal free resolution of the form
	\begin{equation}\label{eqn:bigradedRes}
	0 \lto M \lto \bigoplus_{i,j\in \bbz} \sus^iA(-j)^{\beta_{0,i,j}} \lto \bigoplus_{i,j\in \bbz} \sus^iA(-j)^{\beta_{1,i,j}} \lto \cdots,
	\end{equation}
where the integers $\beta^A_{pij}(M) =\beta_{pij}$ are uniquely determined by $M$ and are called the \textit{trigraded Betti numbers} of $M$. They are all finite, and they can be measured via the dimensions of certain $\Tor$-spaces:
	\[\beta^A_{pij}(M) = \dim_k{ \Tor^A_p{(k,M)}_{ij}}.
	\]
We define the \textit{trigraded Poincar\'e series} of $M$ to be the formal series
	\[\upP^A_M(r,s,t) = \sum_{p,i,j} \beta^A_{pij}(M) r^p s^i t^j.
	\]

For later applications, we note the following proposition which shows that the ``shape'' of the augmentation ideal $\mfm$ is reflected in the ``shape'' of the trigraded Betti numbers of $k$.

\begin{prop}\label{prop:shape}
Let $A$ be a connected bigraded $k$-algebra. If
	\[ \mfm_{ij} \neq 0 \quad \Rightarrow \quad i,j-i>0,
	\]
then
	\[\beta^A_{pij}(k) \neq 0 \quad \Rightarrow \quad i,j-i\geq p.
	\]
\end{prop}

\subsection{Strand totalizations} \label{subsec:strandTotal}

Let $V = \bigoplus_{i,j\in \bbz} V_{ij}$ be a bigraded $k$-vector space. We define an exact functor $V \mapsto V'$ from the category of $\bbz^2$- to $\bbz$-graded vector spaces via the formulas
	\[V' = \bigoplus_{n\in \bbz} V'_n, \quad V'_n = \bigoplus_{j-i=n}V_{ij}.
	\]
The $\bbz$-graded vector space $V'$ is called the \textit{strand totalization} of $V$, and we define the \textit{strand degree} of an element $v\in V_{ij}$ to be the difference $j-i$. The functor $(-)'$ carries a bigraded $k$-algebra $A$ to the graded $k$-algebra $A'$, and it sends a bigraded $A$-module $M$ to the graded $A'$-module $M'$. Its action on shifts is explained in the formula
	\[\left(\sus^pM(q)\right)' = M'(p+q).
	\]

We will say that a connected bigraded $k$-algebra $A$ \textit{lives in positive strands} if its strand totalization $A'$ is connected (as a $\bbz$-graded algebra). If $M$ is a finite bigraded module over such an algebra and $M'$ is a finite $A'$-module, then applying $(-)'$ to the bigraded minimal free resolution \eqref{eqn:bigradedRes} produces the exact sequence
	\[0 \lto M' \lto \bigoplus_{i,j\in \bbz} A'(i-j)^{\beta_{0,i,j}} \lto \bigoplus_{i,j\in \bbz} A'(i-j)^{\beta_{1,i,j}} \lto \cdots,
	\]
which is a graded minimal free resolution of $M'$ over $A'$. This yields the following formula relating the bigraded Betti numbers of $M'$ and the trigraded Betti numbers of $M$:
	\begin{equation}\label{eqn:biTri}
	\beta^{A'}_{pq}(M') = \sum_{j-i = q} \beta_{pij}^A(M).
	\end{equation}
This leads us to

\begin{defn}
Let $A$ be a connected bigraded $k$-algebra. We will say $A$ is \textit{strand-Koszul} if it lives in positive strands and $A'$ is Koszul as a graded algebra.
\end{defn}

Then, in view of equation \eqref{eqn:biTri}, we have

\begin{prop}\label{prop:strandKoszul}
Let $A$ be a connected  bigraded $k$-algebra which lives in positive strands. The algebra $A$ is strand-Koszul if and only if $\beta^A_{pij}(k) \neq 0$ implies $p=j-i$.
\end{prop}

\subsection{DG algebras and modules}\label{subsec:dg}
A \textit{DG $k$-algebra} is a bigraded $k$-algebra $A$ equipped with a $k$-linear endomorphism $\bd^A : A \to A$ of bidegree $(-1,0)$, called the \textit{differential}, such that $\bd^A \circ \bd^A=0$, $\bd^A(1)=0$, and which satisfies the \textit{Leibniz rule}:
	\[\bd(ab) = \bd(a)b +(-1)^{|a|} a \bd (b) \quad \tn{for all $a,b\in A$.}
	\]
We write $A^\natural$ for the underlying bigraded $k$-algebra obtained by forgetting the differential of $A$. If no confusion will arise, we will omit the superscript from $\bd^A$ and write $\bd$ in its place.

The Leibniz rule and the condition $\bd^A(1)=0$ imply that the space $\upZ(A)=\ker{\bd^A}$ of cycles is a bigraded $k$-subalgebra of $A$. The space $\upB(A) = \im{\bd^A}$ of boundaries is an ideal of $\upZ(A)$, and hence the homology $\upH(A) = \upZ(A)/\upB(A)$ is a bigraded $k$-algebra.

Letting $A$ and $B$ be two DG algebras, a function $\varphi:A \to B$ is called a \textit{morphism} of DG algebras if it induces a morphism $A^\natural \to B^\natural$ of the underlying bigraded algebras and if $\bd^B \circ \varphi = \varphi \circ \bd^A$. The morphism $\varphi$ is called a \textit{quasi-isomorphism} if the induced map $\upH(\varphi) : \upH(A) \to \upH(B)$ of bigraded algebras is an isomorphism.

Let $A$ be a DG $k$-algebra.

A \textit{(left) DG $A$-module} is a bigraded left $A^\natural$-module $M$ equipped with a $k$-linear endomorphism $\bd^M: M \to M$ of bidegree $(-1,0)$, called a \textit{differential}, such that $\bd^M \circ \bd^M=0$ and which satisfies the \textit{Leibniz rule}:
	\[\bd^M(am) = \bd^A(a)m +(-1)^{|a|} a \bd^M (m) \quad \tn{for all $a\in A$, $m\in M$.}
	\]
We write $M^\natural$ for the underlying bigraded $A^\natural$-module. Right DG $A$-modules are defined analogously. The homology $\upH(M) = \ker{\bd^M}/\im{\bd^M}$ is a module over $\upH(A)$; the DG module $M$ will be called \textit{homologically finite} if the homology $\upH(M)$ is finite as a bigraded $\upH(A)$-module.

The field $k$ can view viewed as a DG $k$-algebra concentrated in bidegree $(0,0)$ with trivial differential.

If $M$ is a right DG $A$-module and $N$ a left DG $A$-module, their \textit{tensor product}, denoted $M \otimes_A N$, has underlying bigraded $k$-module
	\[(M \otimes_A N)^\natural = M ^\natural \otimes_{A^\natural} N^\natural
	\]
and differential defined as $\bd( m \otimes n) = \bd m \otimes n + (-1)^{|m|} m \otimes \bd n$. Tensor product yields an additive bifunctor $- \otimes_A -: \sfm(A)^\op \times \sfm(A) \to \sfm(k) $ where $\sfm(A)$ and $\sfm(k)$ are the abelian categories of bigraded DG $A$- and $k$-modules, respectively.

Letting $M$ and $N$ be two DG $A$-modules, a function $\alpha:M \to N$ is called a \textit{morphism} of DG modules if it induces an $A^\natural$-linear map $M^\natural \to N^\natural$ and if $\bd^N \circ \alpha = \alpha \circ \bd^M$. The morphism $\alpha$ is called a \textit{quasi-isomorphism} if the induced $\upH(A)$-linear map $\upH(\alpha) : \upH(M) \to \upH(N)$ is an isomorphism.

We will say that $A$ is \textit{connected} if the bigraded algebra $A^\natural$ is connected and there is an inclusion $\bd A_1 \subseteq \mfm_{A_0}$, where $\mfm_{A_0}$ denotes the augmentation ideal of the connected $\bbz$-graded algebra $A_0$. For such a DG algebra, we define the \textit{augmentation ideal} of $A$ to be
	\[\mfm_A = \mfm_{A_0} \oplus A_1 \oplus A_2 \oplus \cdots.
	\]
If $A$ is connected, the natural map $A\to k$ (called the \textit{augmentation}) is a morphism of DG $k$-algebras and gives $k$ a DG $A$-module structure.

\subsection{Semifree resolutions and the derived tensor product}\label{subsec:semifreeRes}

Let $A$ be a DG $k$-algebra. A DG module $F$ is said to be \textit{semifree} if there exists an exhaustive filtration
	\[0 = F^{-1} \subseteq F^0  \subseteq \cdots \subseteq F^{p-1} \subseteq F^p \subseteq \cdots \subseteq F
	\]
by DG submodules such that for each $p\geq 0$ the $A^\natural$-module $(F^p)^\natural/(F^{p-1})^\natural$ is free on a basis of cycles. Such a filtration is called a \textit{semifree filtration} of $F$. For each DG $A$-module $M$, there exists a quasi-isomorphism $\dev:F \to M$ with $F$ a semifree DG $A$-module, and this property determines $F$ up to homotopy equivalence; see \cite[Proposition 6.6]{FelixHalperinThomas2001} for proofs. Such a quasi-isomorphism is called a \textit{semifree resolution} of $M$, though we shall often refer to just $F$ as a semifree resolution of $M$ and omit mention of the map $\dev$.

After choosing a semifree resolution $F(M)$ for every DG module $M$, the left derived functor of $- \otimes_A -$ is defined to be
	\[M \dotimes_A N = F(M) \otimes_A F(N).
	\]
We set
	\[\Tor^A_i{(M,N)}_j = \upH_{ij}(M \dotimes_AN)
	\]
for each $i$ and $j$. By \cite[Proposition 6.7]{FelixHalperinThomas2001}, we have isomorphisms
	\[\upH(F(M) \otimes_A N) \cong \Tor^A{(M,N)} \cong \upH(M \otimes_A F(N))
	\]
of bigraded $k$-modules.

Let $A$ be a connected DG $k$-algebra and $M$ a homologically finite DG $A$-module. Via a process similar to the familiar one used to construct minimal graded free resolutions over a connected algebra, one can construct a semifree resolution $F$ of $M$ such that $F^\natural$ is a finite $A^\natural$-module and $\bd F \subseteq \mfm_A F$. Such a resolution is called a \textit{minimal semifree resolution} of $M$ and is determined by $M$ up to isomorphism of DG $A$-modules. We have that
	\[\Tor^A{(k,M)} \cong k \otimes_A F
	\]
if $F$ is a minimal semifree resolution of $M$. We define the \textit{$(i,j)$-th bigraded Betti number} of $M$ to be
	\[\beta^A_{ij}(M) = \dim_k{\Tor^A_i{(k,M)}_j},
	\]
and we define the \textit{bigraded Poincar\'e series} of $M$ to be the formal series
	\[\upP^A_M(s,t) = \sum_{i,j} \beta^A_{ij}(M) s^i t^j.
	\]

We now offer a definition of Koszulness which applies to connected DG algebras; it is an analog of the definition given by He and Wu in \cite{HeWu2008}.

\begin{defn}
Let $A$ be a connected DG algebra over a field $k$. We say that $A$ is a \textit{Koszul} DG $k$-algebra if $\beta^A_{ij}(k)\neq 0$ implies $i=j$.
\end{defn}

\subsection{Standard graded commutative algebras and Koszul complexes}\label{subsec:kosComp}

A \textit{standard graded $k$-algebra} is a connected graded algebra $R$ such that there exists an isomorphism $R\cong Q/J$ of $k$-algebras where $J$ is a homogeneous ideal in a polynomial ring $Q = k[X_1,\ldots,X_n]$ generated by variables of degree $1$. We assume that $J$ is contained in $\mfm_Q^2$, and in this case the \textit{embedding dimension} of $R$ is the integer $n$. This integer does not depend on the presentation $R\cong Q/J$, since $n = \dim_k{(\mfm_R/\mfm_R^2)}$.

Let $R$ be a standard graded $k$-algebra of embedding dimension $n$, and choose a minimal generating set $\{x_1,\ldots,x_n\}$ of the ideal $\mfm_R$. We view $R$ as a bigraded algebra concentrated in homological degree $0$. Let $t_1,\ldots,t_n$ be a set of exterior variables of bidegree $(1,1)$, and define the \textit{Koszul complex} of $R$ to be the exterior algebra $K^R = \mathsf{\Lambda}_R(t_1,\ldots,t_n)$. We define a differential on $K^R$ by setting $\bd t_i = x_i$ for each $i=1,\ldots,n$, and extending to all of $K^R$ via the Leibniz rule. Equipped with this differential, the algebra $K^R$ is a connected DG $k$-algebra; it is independent (up to isomorphism of DG algebras) of the choice of minimal generating set of $\mfm_R$.

We define the \textit{Koszul homology algebra} of $R$ to be the homology algebra $H^R = \upH(K^R)$; it is a connected bigraded $k$-algebra which lives in positive strands.

If the algebra $R$ under consideration is understood and no confusion is likely to arise, we will write $K$ and $H$ in place of $K^R$ and $H^R$.

\section{Koszulness of $R$, $K^R$, and $H^R$}\label{sec:theory}

Throughout the rest of this paper we let $R$ denote a standard graded $k$-algebra of embedding dimension $n$. We write $K$ and $H$ for its Koszul complex and Koszul homology algebra.

We set as our first task to show that the graded minimal free resolution of $k$ over $R$ and the minimal semifree resolution of $k$ over $K$ are ``the same.'' We begin by (very) briefly recalling a construction due to Tate; for more details, see the original paper \cite{Tate1957}, or \cite[\S6]{Avramov2010}.

\subsection{Tate resolutions}\label{subsec:tate}

The construction begins with the Koszul complex $K$. If $T_1 = \{t_1,\ldots,t_n\}$ is the set of exterior variables generating $K$, we set the notation
	\[R\ang{T_1} = K.
	\]

Now choose cycles $\zeta_{1},\ldots,\zeta_{c}\in\upZ_1(R\ang{T_1})$ whose homology classes minimally generate the module $\upH_1(R\ang{T_1})$. Letting $T_2=\{t_{n+1},\ldots,t_{n+c}\}$ be a set of divided-power variables with $\bideg(t_{n+i}) = (2,\deg{\zeta_i})$ for each $i=1,\ldots,c$, we set
	\[R\ang{T_1,T_2} = R\ang{T_1} \otimes_R \mathsf{\Gamma}_R(T_2), \quad \bd t_{n+i}^{(j)} = \zeta_it^{(j-1)},
	\]
and extend the differential from $R\ang{T_1}$ to all of $R\ang{T_1,T_2}$ via the Leibniz rule. By construction we have
	\[\upH_1(R\ang{T_1,T_2}) =0.
	\]

One then adjoins to $R\ang{T_1,T_2}$ a set $T_3$ of exterior variables in homological degree $3$ in order to ``kill'' cycles whose homology classes minimally generate $\upH_2(R\ang{T_1,T_2})$: the result is the DG algebra
	\[R\ang{T_1,T_2,T_3} = R\ang{T_1,T_2} \otimes_R \mathsf{\Lambda}_R(T_3)
	\]
with $\upH_2(R\ang{T_1,T_2,T_3})=0$. This process continues, alternating between adjoining sets $T_{2i}$ of divided-power variables (in even homological degree) and sets $T_{2i+1}$ of exterior variables (in odd homological degree) to produce a DG algebra 
	\[R\ang{T}, \quad T = \textstyle\bigcup_{i\geq 1} T_i
	\]
which is a graded free resolution of $k$ over $R$. But even more, Gulliksen \cite{Gulliksen1968} and Schoeller \cite{Schoeller1967} proved that this construction yields the \textit{minimal} graded free resolution.

\subsection{A proof of Theorem A}

Using Tate's construction, the proof of the next proposition goes quickly.

\begin{prop}\label{prop:resRRK}
Let $F$ be a minimal graded free resolution of $k$ over $R$ constructed according to Tate's method; see subsection \ref{subsec:tate}.
\begin{enumerate}
\item The resolution $F$ has a DG $K$-module structure making it the minimal semifree resolution of $k$ over $K$.
\item There are equalities
	\[ \upP^R_k(s,t) = (1+st)^n\upP^{K}_k(s,t) \quad \tn{and} \quad \upH_R(t)(1-t)^n\upP^{K}_k(-1,t) = 1,
	\]
where $\upH_R(t) = \sum_{i\geq 0} (\dim{R_i})t^i$ is the Hilbert series of $R$.
\end{enumerate}
\begin{proof}
We carry over the notation introduced in subsection \ref{subsec:tate}, supposing that $F=R\ang{T}$.

(1): The inclusion $K \subseteq R\ang{T}$ exhibits $K$ as a DG subalgebra of $R\ang{T}$, and thus gives $R\ang{T}$ a DG $K$-module structure. However, we have $R\ang{T} = K \otimes_R R\ang{T_{\geq 2}}$, and since $R\ang{T_{\geq 2}}$ is a free bigraded $R$-module, we conclude that $R\ang{T}$ is a semifree DG $K$-module. Minimality of $R\ang{T}$ over $K$ follows immediately from minimality of $R\ang{T}$  over $R$.

(2): We have that
	\[\Tor^R{(k,k)} \cong k \otimes_R R\ang{T} \cong k \ang{T}
	\]
and
	\[\Tor^{K}{(k,k)} \cong k \otimes_K R\ang{T} \cong k \otimes_K ( K \otimes_R R\ang{T_{\geq 2}}) \cong k \ang{T_{\geq 2}}.
	\]
Thus we have
	\[\Tor^R{(k,k)} \cong \mathsf{\Lambda}_k(T_1) \otimes_k \Tor^K{(k,k)}
	\]
as bigraded vector spaces, from which the first equality follows. The second then follows from the first in view of the well-known equality $\upH_R(t)\upP^R_k(-1,t)=1$; see, for example, \cite[Chapter 2, Proposition 2.1]{PolishchukPositselski2005}.
\end{proof}
\end{prop}

As a corollary, we obtain Theorem A from the introduction:

\begin{cor}\label{cor:RKKoszul}
The algebra $R$ is Koszul if and only if $K$ is Koszul.
\begin{proof}
The algebra $R$ is Koszul if and only if
	\[\upP^R_k(s,t) = \sum_{i\geq 0} \beta^R_{ii}(k)(st)^i,
	\]
while $K$ is Koszul if and only if
	\[\upP^K_k(s,t) = \sum_{i\geq 0} \beta^K_{ii}(k)(st)^i.
	\]
The corollary thus follows from the first equation in Proposition \ref{prop:resRRK}(2).
\end{proof}
\end{cor}

Having addressed the relationship between Koszulness of $R$ and $K$, we bring the homology $H$ into the picture. The main vehicle for transporting information from $H$ to $R$ and $K$ will be a pair of spectral sequences.

\subsection{Eilenberg-Moore spectral sequences and quasi-formality}\label{subsec:EMresSS}

Let $\dev: F\simeq k$ be a semifree resolution of $k$ over $K$; recall that this means $\dev: F \to k$ is a quasi-isomorphism and $F$ is a semifree $K$-module. Using a semifree filtration of $F$ (see subsection \ref{subsec:semifreeRes}), one can prove the existence of a collection $\{V^i\}_{i\geq0}$ of bigraded $k$-spaces such that for every $p\geq 0$ we have
	\[(F^p)^\natural \cong K^\natural \otimes_k \textstyle\bigoplus_{i=0}^p V^i \quad \tn{and} \quad \bd(V^p) \subseteq F^{p-1}.
	\]
Thus the spectral sequence arising from the semifree filtration has
	\[E^1_{pqj} \cong (H \otimes_k \sus^{-p}V^p)_{qj}.
	\]
The differentials
	\[d^1_{pqj}: E^1_{pqj} \to E^1_{p-1,q,j}
	\]
and the augmentation $\dev$ assemble together to yield a chain complex
	\[ 0 \lto k = \upH(k) \xleftarrow{\upH(\dev)} H \otimes_k V^0 \xleftarrow{d_1} H \otimes_k \sus^{-1}V^1 \xleftarrow{d_2} H \otimes_k \sus^{-2}V^2 \lto \cdots 
	\]
of free bigraded $H$-modules; if this chain complex is exact, then the semifree resolution $\dev:F \simeq k$ is called an \textit{Eilenberg-Moore resolution} of $k$ over $K$, and the filtration $\{F^p\}_{p\geq -1}$ is called an \textit{Eilenberg-Moore filtration}. These always exist; see, for example, \cite[Proposition 20.11]{FelixHalperinThomas2001}.

From the existence of an Eilenberg-Moore resolution of $k$ over $K$, one can derive the corresponding \textit{Eilenberg-Moore spectral sequence}, which reads as
	\[E^2_{pqj} = \Tor^H_p{(k,k)}_{qj} \Rightarrow \Tor^K_{p+q}{(k,k)}_j
	\]
and with differentials acting as
	\[d^r_{pqj} : E^r_{pqj} \to E^r_{p-r,q+r-1,j}.
	\]
See \cite[\S20]{FelixHalperinThomas2001}, or, alternatively, one can derive this spectral sequence by filtering the bar construction of $K$ by the ``number of bars.'' Adopting Positselski's terminology in \cite{Positselski2017}, we shall say the DG algebra $K$ is \textit{quasi-formal} if the Eilenberg-Moore spectral sequence degenerates on the second page, i.e., all differentials $d^{\geq 2}_{pqj}$ vanish. The terminology is justified by noting that formal DG algebras (i.e., ones which are connected to their homology algebras via a string of quasi-isomorphisms) are quasi-formal; see \cite[Proposition 2.1]{Positselski2017}.

\subsection{Avramov spectral sequences}\label{subsec:AvramovSS}

In addition to the Eilenberg-Moore spectral sequence which links $H$ and $K$, we will use a second spectral sequence which goes directly from $H$ to $R$. To describe it, we first recall that $R$ can be presented as a quotient $Q/J$ where $Q$ is a polynomial ring of embedding dimension $n$. Then, using the classical ``external homology product'' (see \cite{MacLane1963}), one can give $\Tor^Q{(k,k)}$ an $H$-module structure; see also \cite[\S2.3]{Avramov2010}. Our second spectral sequence then reads
	\[E^2_{pqj} = \Tor^H_p{(k,\Tor^Q{(k,k)})}_{qj} \Rightarrow \Tor^R_{p+q}{(k,k)}_j
	\]
and has differentials acting as
	\[d^r_{pqj} : E^r_{pqj} \to E^r_{p-r,q+r-1,j}.
	\]
We call this the \textit{Avramov spectral sequence}, though we note that this name is not standard. For example, Roos \cite{Roos2016} refers to the Eilenberg-Moore spectral sequence above as the Avramov spectral sequence. See \cite{Avramov1981a} for Avramov's original derivation of the spectral sequence, and also Iyengar's alternate derivation in \cite{Iyengar1997}.

Since $Q$ is a polynomial ring, the bigraded $k$-space $\Tor^Q{(k,k)}$ is ``concentrated along the diagonal,'' i.e., $\Tor^Q{(k,k)} = \bigoplus_{i=0}^n \Tor^Q_i{(k,k)}_i$. In view of this fact, the proof of the following result is straightforward.

\begin{prop}\label{prop:strange}
There is an isomorphism
	\[\Tor^Q{(k,k)} \cong \textstyle\bigoplus_{i=0}^n \sus^i k(-i)^{\binom{n}{i}}
	\]
of left bigraded $H$-modules. In particular, there are $k$-linear isomorphisms
	\[\Tor^H_p{(k,\Tor^Q{(k,k)})}_{qj} \cong \textstyle\bigoplus_{i=0}^n \Tor^H_p{(k,k)}_{q-i,j-i}^{\binom{n}{i}}
	\]
for all $p,q,j$ and an equality
	\[\upP^H_{\Tor^Q{(k,k)}}(s,s,t) = (1+st)^n \upP^H_k(s,s,t).
	\]
\end{prop}

\subsection{A proof of Theorem B}

A convergent spectral sequence $E^2_{pqj} \Rightarrow E$ of trigraded vector spaces generates inequalities of sums of dimensions. Indeed, for all $m,j$ we have
\begin{equation}\label{eqn:inequalityDim}
\sum_{p+q=m} \dim_k{E^2_{pqj}} \geq \sum_{p+q=m} \dim_k{E^\infty_{pqj}} = \dim_k{E_{mj}}.
\end{equation}
The following proposition follows.

\begin{prop}\leavevmode\label{prop:inequalities}
\begin{enumerate}
\item The Eilenberg-Moore spectral sequence yields a coefficient-wise inequality
	\[\upP^{K}_k(s,t) \leq \upP^H_k(s,s,t),
	\]
and equality holds if and only if $K$ is quasi-formal.
\item The Avramov spectral sequence yields a coefficient-wise inequality
	\[\upP^R_k(s,t) \leq \upP^H_{\Tor^Q{(k,k)}}(s,s,t),
	\]
and equality holds if and only if the spectral sequence degenerates on the second page.
\end{enumerate}
\end{prop}

As a corollary, we get that the Eilenberg-Moore and Avramov spectral sequences degenerate on the second page simultaneously:

\begin{cor}\label{cor:simulDegen}
The algebra $K$ is quasi-formal if and only if the Avramov spectral sequence degenerates on the second page.
\begin{proof}
Propositions \ref{prop:resRRK}, \ref{prop:strange}, and \ref{prop:inequalities} give
	\[(1+st)^n \upP^K_k(s,t) = \upP^R_k(s,t) \leq \upP^H_{\Tor^Q{(k,k)}}(s,s,t) = (1+st)^n \upP^H_k(s,s,t).
	\]
Now apply Proposition \ref{prop:inequalities} one more time.
\end{proof}
\end{cor}

We now arrive at the central theoretical result of the paper, Theorem B from the introduction, renamed here as

\begin{thm}\label{thm:main}
The following statements are equivalent.
\begin{enumerate}
\item The algebra $H$ is strand-Koszul.
\item The algebra $K$ is Koszul and quasi-formal.
\item The algebra $R$ is Koszul and $K$ is quasi-formal.
\item The algebra $R$ is Koszul and the minimal graded free resolution of $k$ over $R$ has an Eilenberg-Moore filtration.
\end{enumerate}
Furthermore, if one of the above statements is true (and hence all are), then there are equalities
	\[\upP^{K}_k(s,t) = \upP^H_k(s,s,t) \quad \tn{and} \quad \upP^R_k(s,t) = (1+st)^n\upP^H_k(s,s,t).
	\]
\begin{proof}
(1) $\Rightarrow$ (2): We shall first prove that $K$ is quasi-formal. To do so, we observe that if the differential
	\[d^2_{pqj} : E^2_{pqj} \to E^2_{p-2,q+1,j}
	\]
in the Eilenberg-Moore spectral sequence were nonzero, then necessarily
	\[\Tor^H_p{(k,k)}_{qj} = E^2_{pqj} \neq 0 \quad \tn{and} \quad \Tor^H_{p-2}{(k,k)}_{q+1,j} = E^2_{p-2,q+1,j} \neq 0.
	\]
But $H$ is strand-Koszul, and thus Proposition \ref{prop:strandKoszul} implies both $p=j-q$ and $p-2 = j-(q+1)$. This is absurd, which means that we must have $d^2_{pqj}=0$ and therefore $E^2 = E^\infty$.

Now, the coefficient-wise inequality in Proposition \ref{prop:inequalities}(1) turns into an equality, and it yields equalities
	\[\beta^K_{ij}(k) = \sum_{p+q=i} \beta^H_{pqj}(k) \tag{$\ast$}
	\]
for each $i$ and $j$. Thus if $\beta^K_{ij}(k)\neq 0$, then there is a pair $p,q$ with $p+q=i$ and $\beta^H_{pqj}(k) \neq 0$. Proposition \ref{prop:strandKoszul} then gives $p = j-q$, and hence $i=j$. This proves $K$ is Koszul.

(2) $\Rightarrow$ (1): From Proposition \ref{prop:inequalities} and degeneracy of the spectral sequence, we still have the equalities ($\ast$). Thus if $\beta^H_{pqj}(k) \neq 0$ for some $p,q,j$, then $\beta^K_{p+q,j}(k)\neq 0$. Since $K$ is Koszul, this means that $p=j-q$, and by Proposition \ref{prop:strandKoszul}, this means $H$ is strand-Koszul.

(2) $\Leftrightarrow$ (3): Immediate from Corollaries \ref{cor:RKKoszul} and \ref{cor:simulDegen}.

(1) $\Rightarrow$ (4): One can begin from \textit{any} bigraded free resolution
	\[0 \lto k \lto H \otimes_k W^0 \lto H \otimes_k W^1 \lto \cdots
	\]
of $k$ over $H$ (where the $W^p$'s are bigraded $k$-spaces) and produce an Eilenberg-Moore resolution $F$ of $k$ over $K$; see \cite[Proposition 20.11]{FelixHalperinThomas2001}. Indeed, one sets
	\[W = \textstyle\bigoplus_{p=0}^\infty \sus^pW^p \quad \tn{and} \quad F^\natural = K^\natural \otimes_k W,
	\]
and defines the differential on $F$ by induction. If $H$ is strand-Koszul, then we may as well begin with the minimal bigraded free resolution of $k$ over $H$, in which case $(\sus^pW^p)_{p+i,j}\neq 0$ implies $p+i=j$ for all $p,i,j$ (see Proposition \ref{prop:strandKoszul}). This means that the basis of $F^\natural$ is concentrated along the diagonal, and thus $F$ is the minimal semifree resolution of $k$ over $K$.

(4) $\Rightarrow$ (1): Let $F$ be the minimal semifree resolution of $k$ over $K$ and let $\{F^p\}_{p\geq -1}$ be an Eilenberg-Moore filtration. Since $K$ is Koszul, in the notation of subsection \ref{subsec:EMresSS} the bigraded vector spaces $\{V^i\}_{i\geq 0}$ are concentrated along the diagonal. But then the resulting free resolution
	\[0 \lto k \lto H \otimes_k V^0 \lto H \otimes_k \sus^{-1}V^1 \lto H \otimes_k \sus^{-2}V^2 \lto \cdots 
	\]
of $k$ over $K$ has the property that $(\sus^{-p}V^p)_{ij}\neq 0$ implies $p= j-i$ for all $p,i,j$. Hence $H$ is strand-Koszul, by Proposition \ref{prop:strandKoszul}.
\end{proof}
\end{thm}

\begin{example}\label{ex:CI}
One class of algebras for which we can explicitly write down Eilenberg-Moore filtrations is the class of quadratic complete intersections. If $R$ is such an algebra, then $H$ is strand-Koszul; this follows, for example, from a (graded version of a) result of Tate \cite{Tate1957}, which states that $H$ is a bigraded exterior algebra generated by $H_{1,2}$; hence the stand totalization $H'$ (recall the definition in section \ref{subsec:strandTotal}) is an exterior algebra generated in degree $1$, which is well-known to be Koszul.

The details are as follows. Assume $R$ is a quadratic complete intersection of embedding dimension $n$ and codimension $c$; this means that
	\[R \cong \frac{k[X_1,\ldots,X_n]}{(G_1,\ldots,G_c)}
	\]
where $G_1,\ldots,G_c$ is a regular sequence of degree-$2$ homogeneous polynomials. Write $x_i$ for the image of $X_i$ in $R$, so that the Koszul complex is $K = \mathsf{\Lambda}_R(t_1,\ldots,t_n)$ with $\bd t_i = x_i$ for each $i$. If we write
	\[G_h = \sum_{1 \leq i \leq j \leq n} \lambda_{hij} X_i X_j, \quad \lambda_{hij} \in k
	\]
for each $h=1,\ldots,c$, then each $z_h = \sum_{1\leq i \leq j \leq n} \lambda_{hij} x_i t_j$ is a cycle in $K$. In fact, if we write $\zeta_h$ for the homology class of $z_h$, then Tate's result mentioned above shows that
	\[H = \mathsf{\Lambda}_k(\zeta_1,\ldots,\zeta_c).
	\]

As always, a minimal graded free resolution $F$ of $k$ over $R$ can be constructed according to Tate's method; see subsection \ref{subsec:tate}. However, since $R$ is a complete intersection, Tate's theory shows that the construction stops after adjoining divided-power variables in homological and internal degree $2$; the result is that
	\[F = K \otimes_R \mathsf{\Gamma}_R(s_{1},\ldots,s_{c}), \quad \bd s_h = z_h,
	\]
where $s_{1},\ldots,s_{c}$ are divided-power variables of bidegree $(2,2)$.

Now, for each $p\geq 0$, set
	\[F^p = \sum_{i_1 + \cdots + i_c\leq p} K s_1^{(i_1)} \cdots s_c^{(i_c)} \quad \tn{and} \quad V^p = \sum_{i_1 + \cdots + i_c= p} k s_1^{(i_1)} \cdots s_c^{(i_c)}.
	\]
One can check that $\{F^p\}_{p\geq 0}$ is an increasing filtration of $F$ by DG $K$-submodules which is semifree. In the resulting spectral sequence the differential $d^1_{p}: E^1_{p} \to E^1_{p-1}$
can be identified with the mapping 
	\[H \otimes_k \sus^{-p}V^p \to H \otimes_k \sus^{-p+1}V^{p-1}
	\]
defined on a basis element by
	\[1 \otimes \sus^{-p} s_1^{(i_1)} \cdots s_c^{(i_c)} \mapsto \sum_{j=1}^c \zeta_j \otimes \sus^{-p+1} s_1^{(i_1)} \cdots s_j^{(i_j-1)}\cdots s_c^{(i_c)}
	\]
where we set $s^{(i_j-1)}= 0$ if $i_j-1<0$. Thus
	\[ 0 \lto H \otimes_k V^0 \xleftarrow{d_1} H \otimes_k \sus^{-1}V^1 \xleftarrow{d_2} H \otimes_k \sus^{-2}V^2 \lto \cdots 
	\]
is the minimal bigraded free resolution of $k$ over $H$; one can verify this, for example, by using Koszul duality theory (see \cite[Chapter 2, \S3]{PolishchukPositselski2005}). We conclude that the filtration $\{F^p\}$ is an Eilenberg-Moore filtration.
\end{example}

We finish this section by noting that even when the Avramov spectral sequence does not degenerate predictably, it is still possible to obtain relations from it between the low-degree Betti numbers of $k$ over $R$ and $k$ over $H$. These relations are presented in

\begin{thm}
The following equalities hold:
\begin{align*}
\beta^R_{2,2}(k) &= \binom{n}{2} +\beta^H_{1,1,2}(k),\\
\beta^R_{2,j}(k) &= \beta^H_{1,1,j}(k)\quad \tn{for $j>2$,}\\
\beta^R_{3,j}(k) &= \beta^H_{1,2,j}(k) + n \beta^H_{1,1,j-1}(k) \quad \tn{for $j>3$,} \\
\beta^R_{4,j}(k) &= \beta^H_{1,3,j}(k) + n\beta^H_{1,2,j-1}(k) + \binom{n}{2}\beta^H_{1,1,j-2}(k) + \beta^H_{2,2,j}(k) \quad \tn{for $j>4$.}
\end{align*}
\begin{proof}
Combining the Avramov spectral sequence with Proposition \ref{prop:strange} produces a spectral sequence
	\[E^2_{pqj} = \textstyle\bigoplus_{i=0}^n \Tor^H_p{(k,k)}_{q-i,j-i}^{\binom{n}{i}} \Rightarrow \Tor^R_{p+q}{(k,k)}_j
	\]
with differentials on the second page acting on tridegrees as
	\[d^2_{pqj} : E^2_{pqj} \to E^2_{p-2,q+1,j}.
	\]
The algebra $H$ satisfies the hypotheses of Proposition \ref{prop:shape}; thus if $E^2_{pqj}\neq 0$, then there is an $i$ with $0 \leq i \leq n$ such that $\beta^H_{p,q-i,j-i}(k)\neq 0$. Hence $q-i\geq p$, so that necessarily $q \geq p$. Thus from \eqref{eqn:inequalityDim} we have
	\[\beta^R_{mj}(k) = \sum_{\substack{p+q=m \\ q\geq p}} \dim_k{E^\infty_{pqj}}
	\]
for all $m$ and $j$.

In particular, we have
	\[\beta^R_{2,j}(k) = \dim_k{E^\infty_{0,2,j}} + \dim_k{E^\infty_{1,1,j}}.
	\]
However, according to the description of the differentials $d^2_{pqj}$ given above, we must have $E^\infty_{0,2,j} = E^2_{0,2,j}$ and $E^\infty_{1,1,j} = E^2_{1,1,j}$. Thus
\begin{align*}
\beta^R_{2,j}(k) &= \sum_{i=0}^n \binom{n}{i} \beta^H_{0,2-i,j-i}(k) + \sum_{i=0}^n \binom{n}{i} \beta^H_{1,1-i,j-i}(k) \\
&= \binom{n}{2}\beta^H_{0,0,j-2}(k) + \beta^H_{1,1,j}(k),
\end{align*}
from which the equations for $\beta^R_{2,j}(k)$ follow.

Similarly, we have
\begin{align*}
\beta^R_{3,j}(k) &= \dim_k{E^\infty_{0,3,j}} + \dim_k{E^\infty_{1,2,j}} \\
&= \dim_k{E^\infty_{0,3,j}} + \dim_k{E^2_{1,2,j}} \\
&= \dim_k{E^\infty_{0,3,j}} + \beta^H_{1,2,j}(k) + n \beta^H_{1,1,j-1}(k).
\end{align*}
But
	\[E^2_{0,3,j} = \bigoplus_{i=0}^n \Tor^H_0{(k,k)}^{\binom{n}{i}}_{3-i,j-i} = \Tor^H_0{(k,k)}^{\binom{n}{3}}_{0,j-3},
	\]
so that $E^\infty_{0,3,j} = 0$ if $j>3$. The equations for $\beta^R_{3,j}(k)$ follow.

Finally, we have
\begin{align*}
\beta^R_{4,j}(k) &= \dim_k{E^\infty_{0,4,j}} + \dim_k{E^\infty_{1,3,j}} + \dim_k{E^\infty_{2,2,j}} \\
&= \dim_k{E^\infty_{0,4,j}} + \dim_k{E^2_{1,3,j}} + \dim_k{\ker{d^2_{2,2,j}}}.
\end{align*}
One checks just as above that $E^\infty_{0,4,j} = 0$ and $\ker{d^2_{2,2,j}} = E^2_{2,2,j}$ when $j>4$; the equations for $\beta^R_{4,j}(k)$ then follow.
\end{proof}
\end{thm}

\begin{remark}
The first two equations for $\beta^R_{2,j}(k)$ refine well-known existing relations; see, for example, \cite[Theorem 2.3.2]{BrunsHerzog1998}.
\end{remark}

\section{Classes of Koszul homology algebras which are strand-Koszul}

In this section we describe several classes of standard graded algebras that possess Koszul homology algebras which are strand-Koszul. Our goal is to systematically work through the proof of Theorem C from the introduction, noting that part (3) of that theorem has already been established in Example \ref{ex:CI}.

Throughout this section $R$ denotes a standard graded algebra written as a quotient of a polynomial ring $Q$; see subsection \ref{subsec:kosComp}.

\subsection{Small embedding codepth and Golod algebras}

We define the \textit{depth} of $R$, denoted $\depth{R}$, to be the maximal length of a (homogeneous) regular sequence in the augmentation ideal of $R$.

\begin{prop}
Let $R$ be a standard graded algebra of embedding dimension $n$ and let $H$ be its Koszul homology algebra. If $n-\depth{R}\leq 3$, then the algebra $R$ is Koszul if and only if $H$ is strand-Koszul.
\begin{proof}
Avramov proved in \cite[Theorem 5.9]{Avramov1974} that the Koszul complex $K$ is quasi-formal when $n-\depth{R} \leq 3$. Now apply Theorem \ref{thm:main}.
\end{proof}
\end{prop}

By using bar constructions, Iyengar \cite{Iyengar1997} gave an alternate construction of the Avramov spectral sequence (see subsection \ref{subsec:AvramovSS}) and described its first page as
	\[E^1_{pqj} = \left( \cl{H} ^{\otimes p} \otimes_k \Tor^Q{(k,k)} \right)_{qj},
	\]
where $\cl{H}$ is the cokernel of the unit map $k \to H$. We can thus add a series to the inequality in Proposition \ref{prop:inequalities}(2) to produce
	\begin{equation}\label{eqn:iyengar}
	\upP^R_k(s,t) \leq \upP^H_{\Tor^Q{(k,k)}}(s,s,t) \leq \frac{(1+st)^n}{1-t(\upP^Q_R(s,t)-1)}.
	\end{equation}
If the series at the two ends are equal, then $R$ is called a \textit{Golod algebra}. These rings have been intensively studied in the homology theory of commutative local rings (see \cite{Avramov2010} for an overview).

\begin{prop}
Let $R$ be a standard graded algebra and let $H$ be its Koszul homology algebra. If $R$ is Koszul and Golod, then $H$ is strand-Koszul.
\begin{proof}
If $R$ is a Golod algebra, then the inequalities in \eqref{eqn:iyengar} are equalities, and by Proposition \ref{prop:inequalities} and Corollary \ref{cor:simulDegen} this means the Koszul complex $K$ is quasi-formal. If $R$ is Koszul as well, we can apply Theorem \ref{thm:main} to obtain the desired conclusion.
\end{proof}
\end{prop}

In the rest of this section, the technique we use for establishing strand-Koszulness comes from the theory of noncommutative Gr\"obner bases. Our main reference is \cite{Bergman1978}, but we will collect here the main definitions and results for the reader's convenience.

\subsection{Noncommutative Gr\"obner bases}

Let $k$ be a field and write
	\[T = k \ang{\zeta_1,\ldots,\zeta_n}
	\]
for the polynomial algebra generated by \textit{noncommuting} variables $\zeta_1,\ldots,\zeta_n$ of degree $1$. We consider the degree-lexicographic ordering on the monomials of $T$ with the variables ordered
	\[\zeta_1 < \cdots < \zeta_n.
	\]
Let $G = \{\gamma_1,\ldots,\gamma_t\}$ be a collection of elements of $T$ which generate an ideal $I$. Multiplying through by scalars if necessary, for each $i=1,\ldots,t$ we can write $\gamma_i = \mu_i +\alpha_i$ where $\mu_i$ is the leading monomial of $\gamma_i$ (with coefficient $1$) and $\alpha_i$ is a $k$-linear combination of monomials which are all $<\mu_i$. A nonzero monomial $\mu\in T$ is said to be \textit{reduced} (with respect to $G$) if it does not contain any of the $\mu_i$'s as a sub-monomial; we then say $G$ forms a \textit{Gr\"obner basis} of $I$ if the set of images in $T/I$ of the reduced monomials form a $k$-basis for the quotient algebra. These images always linearly span the quotient algebra, so to verify that a set of generators of $I$ is a Gr\"obner basis we need only check linear independence.

The relevance of Gr\"obner bases is explained by the following fact: If a \textit{quadratic} Gr\"obner basis of $I$ exists, then the quotient $T/I$ is necessarily Koszul (in commutative algebra such algebras are called \textit{$G$-quadratic}). Indeed, this follows from a filtration argument; see, for example, \cite[Chapter 4, Theorem 7.1]{PolishchukPositselski2005}.

\subsection{Short Gorenstein algebras}
Suppose $R$ is artinian Gorenstein of socle degree $2$ and let $H$ be its Koszul homology algebra. We then have that $H_{n,n+2} \cong k$ and $H_{n,j} = 0$ for all $j\neq n+2$. But in fact, more is true: Avramov and Golod proved in \cite{AvramovGolod1971} that the Gorenstein condition implies $H$ is a Poincar\'e algebra, which means that for all $i,j$ there are $k$-linear isomorphisms
	\[H_{ij} \to \Hom_k{(H_{n-i,n+2-j},H_{n,n+2})}
	\]
given by $h\mapsto \lambda_h$ where $\lambda_h$ is (left) multiplication by $h$. In particular, the Betti table of $R$ over $Q$ looks like

\begin{center}
\begin{tabular}{C | C  C  C  C  C C }
& 0 & 1 & 2 & \cdots & n-1 & n \\ \hline
0 & 1 & - & - & \cdots & - & - \\ 
1 & - & b_1 & b_2 & \cdots & b_{n-1} & - \\ 
2 & - & - & - & \cdots & - & 1
\end{tabular}
\end{center}

\noindent where $b_i = b_{n-i}$ for each $i=1,2,\ldots, \lfloor n/2 \rfloor$.

\begin{thm}\label{thm:shortGor}
Let $R$ be an artinian Gorenstein standard graded $k$-algebra of socle degree $2$ and embedding dimension $n$, and let $H$ be the Koszul homology algebra of $R$. If $k$ does not have characteristic $2$, or if $n$ is odd, then $H$ is strand-Koszul.
\begin{proof}
Suppose first that $n$ is odd, set $c = \lfloor n/2 \rfloor$, and consider the multiplication maps
	\[\mu_i:H_{i,i+1} \otimes_k H_{n-i,n-i+1} \to H_{n,n+2}
	\]
for $i=1,2,\ldots,c$. As can be seen from the Betti table above, these are the only nonzero multiplication maps besides the ones that involve $H_{0,0}$. Since $H_{n,n+2}$ is $1$-dimensional, fixing a nonzero element $\sigma \in H_{n,n+2}$ allows us to identify $H_{n,n+2}$ with $k$. Then, by the Poincar\'e condition on $H$, for each $i=1,2,\ldots,c$ there are ordered bases $\zeta_{i,1},\ldots,\zeta_{i,b_i}$ of $H_{i,i+1}$ and $\eta_{n-i,1},\ldots,\eta_{n-i,b_i}$ of $H_{n-i,n-i+1}$ such that
	\begin{equation}\label{eqn:rel1}
	\zeta_{i,j} \eta_{n-i,\ell} = \begin{cases}
	0 & : j\neq \ell, \\
	\sigma & : j = \ell,
	\end{cases}
	\end{equation}
for all $1\leq j,\ell \leq b_i$.

We set $S = \bigcup_{i=1}^{c} \{\zeta_{i,1},\ldots,\zeta_{i,b_i},\eta_{n-i,1},\ldots,\eta_{n-i,b_i}\}$ and we consider the noncommutative polynomial ring $k\ang{S}$ generated by $S$. There is an evident surjective $k$-algebra map
	\[\varphi: k\ang{S} \surj H
	\]
which sends the elements of $S$ to themselves and which preserves bidegrees (and hence also strand degrees). Therefore, if we set $I = \ker{\varphi}$ and write $H'$ for the strand totalization of $H$, then we have an isomorphism $k\ang{S}/I \cong H'$ of graded algebras. Thus to prove $H$ is strand-Koszul it will suffice to show that $I$ has a quadratic (with respect to strand degree) Gr\"obner basis.

Let $G$ be the subset of $k\ang{S}$ consisting of the following four types of elements:
\begin{enumerate}
\item all monomials of strand degree $2$, except for monomials of the form $\zeta_{i,j}\eta_{n-i,j}$ and $\eta_{n-i,j}\zeta_{i,j}$ for $i=1,\ldots,c$ and $j = 1,\ldots,b_i$,
\item all commutators $\eta_{n-i,j}\zeta_{i,j} -(-1)^{i(n-i)} \zeta_{i,j}\eta_{n-i,j}$ for $i=1,\ldots,c$ and $j = 1,\ldots,b_i$,
\item all elements of the form $\zeta_{i,j}\eta_{n-i,j} - \zeta_{1,1}\eta_{n-1,1}$ for $i=2,\ldots,c$ and $j=1,\ldots,b_i$, and
\item all elements of the form $\zeta_{1,j}\eta_{n-1,j} - \zeta_{1,1}\eta_{n-1,1}$ for $j=2,\ldots,b_i$.
\end{enumerate}
If we use the degree-lexicographic ordering on the monomials in $k\ang{S}$ with
\begin{multline*}
\zeta_{1,1} < \cdots < \zeta_{1,b_1} < \cdots < \zeta_{c,1} < \cdots < \zeta_{c,b_c} < \eta_{n-c,1} < \cdots  \\ \cdots < \eta_{n-c,b_c} < \cdots < \eta_{n-1,1} < \cdots < \eta_{n-1,b_1},
\end{multline*}
then the only reduced monomial (with respect to $G$) of strand degree $\geq 2$ is $\zeta_{1,1}\eta_{n-1,1}$. Since this monomial is not in the ideal $(G)$, we conclude that $G$ forms a Gr\"obner basis of $(G)$. Using the Betti diagram of $R$ over $Q$, equation \eqref{eqn:rel1}, and the fact $H$ is a quotient of an exterior algebra, one can show that $G\subseteq I$. Thus we have a surjection
	\[\frac{k\ang{S}}{(G)} \surj H
	\]
induced by $\varphi$. Finally, since we have $\dim{(k\ang{S}/(G))_{ij}} = \dim{H_{ij}}$ for all $i,j$, we have $(G)=I$. Hence $I$ has a quadratic Gr\"obner basis.

If $n$ is even, the same argument as above goes through as long as $k$ does not have characteristic $2$; indeed, the ``middle'' multiplication map
	\[H_{n/2,n/2+1} \otimes H_{n/2,n/2+1} \to H_{n,n+2}
	\]
can then be diagonalized.
\end{proof}
\end{thm}

\subsection{Koszul algebras defined by three relations}

Suppose $R$ is Koszul and that $J$ (the defining ideal of $R$; see subsection \ref{subsec:kosComp}) is minimally generated by three elements. In \cite{BoocherHassanzadehIyengar2017} it is shown that $H$ is generated by its linear strand and that the Betti table of $R$ over $Q$ must be one of the following four:

\medskip
\begin{center}
\begin{tabular}{C | C  C  C  C}
& 0 & 1 & 2 & 3 \\ \hline
0 & 1 & - & - & -\\ 
1 & - & 3 & 2 & -
\end{tabular} \hspace{0.25in}\begin{tabular}{C | C  C  C  C}
& 0 & 1 & 2 & 3 \\ \hline
0 & 1 & - & - & -\\ 
1 & - & 3 & 3 & 1
\end{tabular}
\end{center}

\medskip

\begin{center}
\begin{tabular}{C | C  C  C  C}
& 0 & 1 & 2 & 3 \\ \hline
0 & 1 & - & - & -\\ 
1 & - & 3 & - & -\\ 
2 & - &  - & 3 & - \\
3 & - & - & - & 1
\end{tabular} \hspace{0.25in}\begin{tabular}{C | C  C  C  C}
& 0 & 1 & 2 & 3 \\ \hline
0 & 1 & - & - & -\\ 
1 & - & 3 & 1 & -\\ 
2 & - &  - & 2 & 1  
\end{tabular}
\end{center}

\medskip

\noindent If the Betti table of $R$ is one of the two in the top row, then $H$ is obviously strand-Koszul. The bottom-left Betti table corresponds to an $H$ which is an exterior algebra generated by $H_{1,2}$, and hence $H$ is again strand-Koszul. Thus in the proof of the next theorem we may assume that the Betti table of $R$ is equal to the bottom-right table.

\begin{thm}
Let $R$ be a Koszul standard graded algebra. If the defining ideal of $R$ is minimally generated by three elements, then the Koszul homology algebra of $R$ is strand-Koszul.
\begin{proof}
As we just noted in the paragraph immediately preceding the theorem, the algebra $H$ is generated by its linear strand. Using this fact along with the form of the Betti table of $R$, we can choose elements $\zeta_1,\zeta_2,\zeta_3,\eta,\sigma \in H$ such that each $\zeta_1,\zeta_2,\zeta_3$ has bidegree $(1,2)$, $\eta$ has bidegree $(2,3)$, $\sigma$ has bidegree $(3,5)$, and
	\[\zeta_1\eta  = \sigma, \quad \zeta_2\eta=\zeta_3\eta=0.
	\]
Furthermore, since the $2$-dimensional vector space $H_{2,4}$ is spanned by the three monomials $\zeta_1\zeta_2,\zeta_1\zeta_3,\zeta_2\zeta_3$, there must be $a,b,c\in k$, not all zero, such that
	\[a\zeta_1\zeta_2 + b\zeta_1\zeta_3 + c\zeta_2\zeta_3 =0.
	\]

We have an evident surjective $k$-algebra map
	\[\varphi: k\ang{\zeta_1,\zeta_2,\zeta_3,\eta} \surj H.
	\]
This map preserves bidegrees and strand degrees, so we have a graded $k$-algebra isomorphism $k\ang{\zeta_1,\zeta_2,\zeta_3,\eta}/I \cong H'$ where $I = \ker{\varphi}$. As in the proof of Theorem \ref{thm:shortGor}, we will show $H$ is strand-Koszul by proving $I$ has a quadratic (with respect to strand degree) Gr\"obner basis.
	
Let $G$ be the subset of $k\ang{\zeta_1,\zeta_2,\zeta_3,\eta}$ consisting of the following four types of elements:
\begin{enumerate}
\item all commutators $\eta \zeta_i - \zeta_i \eta$ for $1\leq i \leq 3$,
\item all skew commutators $\zeta_j \zeta_i + \zeta_i \zeta_j$ for $1\leq i < j \leq 3$,
\item all squares $\zeta_1^2$, $\zeta_2^2$, $\zeta_3^2$, $\eta^2$ and $\zeta_2\eta$, $\zeta_3\eta$,
\item a linear combination $a\zeta_1\zeta_2 + b \zeta_1\zeta_3 + c\zeta_2\zeta_3$ where not all $a,b,c,\in k$ are zero.
\end{enumerate}

If $c\neq 0$, then with the degree-lexicographic ordering and
	\[\zeta_1 < \zeta_2 < \zeta_3 < \eta,
	\]
the reduced monomials of degree $\geq 2$ are $\zeta_1 \eta$, $\zeta_1\zeta_2$, and $\zeta_1\zeta_3$. Since the images of these monomials in the quotient $k\ang{\zeta_1,\zeta_2,\zeta_3,\eta}/(G)$ are linearly independent, we conclude that $G$ forms a Gr\"obner basis of $(G)$ and hence the quotient $k\ang{\zeta_1,\zeta_2,\zeta_3,\eta}/(G)$ is Koszul.

If $c=0$ and $a,b\neq 0$, then the reduced monomials of degree $\geq 2$ are $\zeta_1\eta $, $\zeta_1\zeta_2$, and $\zeta_2\zeta_3$. These form a Gr\"obner basis of $(G)$, and the quotient $k\ang{\zeta_1,\zeta_2,\zeta_3,\eta}/(G)$ is again Koszul.

If $c=0$ and one of $a$ or $b$ is also $0$, then $k\ang{\zeta_1,\zeta_2,\zeta_3,\eta}/(G)$ is the quotient of a skew-polynomial algebra by an ideal generated by quadratic monomials; such algebras are Koszul (see \cite[Chapter 4, Theorem 8.1]{PolishchukPositselski2005}).

It is not difficult to prove that $G\subseteq I$, and hence we have a surjection
	\[\frac{k\ang{\zeta_1,\zeta_2,\zeta_3,\eta}}{(G)} \surj H
	\]
induced by $\varphi$. Since the dimensions of $k\ang{\zeta_1,\zeta_2,\zeta_3,\eta}/(G)$ and $H$ are equal in each bidegree, we have that $(G) = I$. Hence $H$ is strand-Koszul.
\end{proof}
\end{thm}

\subsection{Koszul algebras defined by edge ideals of paths}
Finally, we turn our attention to algebras of the form
	\[R = \frac{k[X_1,\ldots,X_n]}{(X_1X_2,X_2X_3,\ldots,X_{n-1}X_n)}, \quad n\geq 3.
	\]
Again letting $H$ denote the Koszul homology algebra of $R$, we shall prove $H$ is strand-Koszul via the same Gr\"obner basis techniques used above and the main input to our proof will again be the dimensions of the homogeneous components of $H$. However, this time our $R$ has extra structure that will prove to be useful; namely, the defining ideal of $R$ is generated by monomials and hence $R$ has an $\bbz^n$-grading which refines its usual $\bbz$-grading. This finer grading induces a $\bbz^{n+1}$-grading on $H$ which, in turn, refines its usual $\bbz^2$-grading. Thus our first goal is to set up some notation and terminology to describe these finer gradings.

\subsubsection{Multidegrees}\stepcounter{equation}

Elements of $\bbz^n$ will be called \textit{multidegrees}. If $\mbv=(v_1,\ldots,v_n)$ is a multidegree, we set
	\[| \mbv | = v_1 + \cdots + v_n
	\]
and we define the \textit{support} of $\mbv$ to be the set
	\[\supp{(\mbv)} = \{ i \mid v_i \neq 0\}.
	\]
We call $\mbv$ \textit{squarefree} if every component $v_i$ is either $0$ or $1$. 

Write
	\[\mbe_1 = (1,0,\ldots,0), \ \mbe_2 = (0,1,0,\ldots,0), \ \ldots \ , \mbe_n = (0,\ldots,0,1).
	\]
Then, for $1 \leq i \leq n$ and $1\leq r \leq n-i+1$, set
	\[\mbp_{i,r} = \mbe_{i} + \cdots + \mbe_{i+r-1} \in \bbz^n.
	\]
Every squarefree multidegree $\mbv\in \bbz^n$ (not equal to $(0,\ldots,0)$) can be decomposed uniquely into a sum of $\mbp_{i,r}$'s with maximal support; for example,
	\[\mbv = (1,1,1,0,1,1,0,0,1,1,1,1) \in \bbz^{12}
	\]
decomposes as $\mbv = \mbp_{1,3} + \mbp_{5,2} + \mbp_{9,4}$. Such a decomposition will be called a \textit{complete decomposition} of $\mbv$.

The polynomial ring $Q = k[X_1,\ldots,X_n]$ is $\bbz^n$-graded by assigning the variable $X_i$ the multidegree $\mbe_i$. The defining ideal $(X_1X_2,X_2X_3,\ldots,X_{n-1}X_n)$ of $R$ is homogeneous with respect to this $\bbz^n$-grading, and hence $R$ inherits a $\bbz^n$-grading from $Q$.

Write $x_i$ for the image of $X_i$ in $R$ and let $K$ be the Koszul complex of $R$ generated by exterior variables $t_1,t_2,\ldots,t_n$. Given multidegrees $\mbv,\mbw\in \bbz^n$, we define monomials
	\[x^\mbv = x_1^{v_1} \cdots x_n^{v_n}\in R \quad \tn{and} \quad t^\mbw = t_1^{w_1} \cdots t_n^{w_n}\in K.
	\]
The monomial $x^\mbv t^\mbw\in K$ is then assigned the \textit{multigraded bidegree} $(| \mbw |, \mbv+\mbw) \in \bbz^{n+1}$, where $|\mbw|$ is the usual homological degree of the monomial and $\mbv+\mbw$ is called the \textit{multigraded internal degree}. The differential on $K$ preserves this latter degree and hence the homology $H$  inherits a $\bbz^{n+1}$-grading from $K$. For each $i\in \bbz$ and $\mbu \in \bbz^n$, we let $K_{i,\mbu}$ and $H_{i,\mbu}$ be the $k$-subspaces of $K$ and $H$, respectively, which are spanned by all homogeneous elements of multigraded bidegree $(i,\mbu)$; for each $j\in \bbz$ we then have decompositions
	\[K_{ij} = \bigoplus_{\substack{\mbu \in \bbz^n \\ |\mbu | = j}} K_{i,\mbu} \quad \tn{and} \quad H_{ij} = \bigoplus_{\substack{\mbu \in \bbz^n \\ |\mbu| = j}} H_{i,\mbu}.
	\]
Hence the multigraded bidegrees on $K$ and $H$ refine the usual bidegrees. The strand degree of a homogeneous element of multigraded bidegree $(i,\mbu)$ is $|\mbu|-i$.

\medskip
We now turn our attention to the result which serves as the basis for our study of $H$:

\begin{prop}[Boocher, D'Al\`\i, Grifo, Monta\~no, Sammartano \cite{BoocherEtAl2015}]\label{prop:BoocherEtAl}\leavevmode
\begin{enumerate}
\item The algebra $H$ is generated by the subspaces $H_{1,2}$ and $H_{2,3}$.
\item Let $\mbu\in \bbz^n$ be a squarefree multidegree with complete decomposition $\mbu = \sum_{p=1}^m \mbp_{i_p,r_p}$. Then
	\[\dim_k{ H_{\ast,\mbu}} = \begin{cases}
	1 & : r_p \not\equiv 1 \modu{3}\tn{ for each $p$}, \\
	0 & : \tn{otherwise}.
	\end{cases}
	\]
In the former case, $\dim_k{H_{i,\mbu}}=1$ exactly when $i = \sum_{p=1}^m \left\lfloor 2r_p/3 \right\rfloor$.
\end{enumerate}
\end{prop}

\begin{remark}\label{rmk:squarefree}
We note that $H_{\ast,\mbu}=0$ if $\mbu$ is not squarefree (see, for example, \cite[Corollary 1.40]{MillerSturmfels2005}), so Proposition \ref{prop:BoocherEtAl} describes the dimensions of \textit{all} nonzero $H_{i,\mbu}$.
\end{remark}

Now, for all $i=1,\ldots,n-1$ and $j=1,\ldots,n-2$, the monomials $x_{i+1}t_i$ and $x_{j+1}t_j t_{j+2}$ are cycles in the Koszul complex $K$, and hence we may define the homology classes
	\[\zeta_{i,i+1} = \cls{(x_{i+1}t_i)} \in H_{1,2} \quad \tn{and} \quad \eta_{j,j+1,j+2} = \cls{(x_{j+1} t_j t_{j+2})} \in H_{2,3}.
	\]
The homology classes $\zeta_{i,i+1}$ and $\eta_{j,j+1,j+2}$ have multigraded bidegrees $(1,\mbp_{i,2})$ and $(2,\mbp_{j,j+1,j+2})$, respectively, and thus they all have strand degree $1$.

\begin{prop}\label{prop:bases}
The sets $\{\zeta_{1,2},\ldots,\zeta_{n-1,n}\}$ and $\{\eta_{1,2,3},\ldots,\eta_{n-2,n-1,n}\}$ form $k$-bases for the spaces $H_{1,2}$ and $H_{2,3}$, respectively.
\begin{proof}
That the first set is a basis of $H_{1,2}$ follows from \cite[Theorem 2.3.2]{BrunsHerzog1998}, for example.

To prove the second set is a basis of $H_{2,3}$, we first note that
	\[H_{2,3} = \bigoplus_{\substack{\mbu \in \bbz^n \\ |\mbu| = 3}} H_{2,\mbu}.
	\]
We claim, in fact, that
	\[H_{2,3} = \bigoplus_{j=1}^{n-2} H_{2,\mbp_{j,3}}. \tag{$\ast$}
	\]
To see this, suppose that $H_{2,\mbu}\neq 0$ and $|\mbu|=3$. By Remark \ref{rmk:squarefree}, we have that $\mbu$ must be squarefree. But then Proposition \ref{prop:BoocherEtAl}(2) shows that $\mbu = \mbp_{j,3}$ for some $j \in \{1,2,\ldots,n-2\}$, which establishes the equation $(\ast)$. But the same proposition shows that each $H_{2,\mbp_{j,3}}$ is $1$-dimensional, and thus to prove that the set $\{\eta_{1,2,3},\ldots,\eta_{n-2,n-1,n}\}$ is a basis of $H_{2,3}$ it will suffice to show that $\eta_{j-2,j-1,j}$ is nonzero (since it has multigraded bidegree $(2,\mbp_{j,3})$).

To do this, first note that $K_{2,\mbp_{j,3}}$ and $K_{3,\mbp_{j,3}}$ are $3$- and $1$-dimensional, respectively, with bases
	\[\{x_{j+1}t_jt_{j+2}, x_j t_{j+1}t_{j+2},x_{j+2}t_jt_{j+1}\} \quad \tn{and} \quad \{ t_j t_{j+1}t_{j+2}\}.
	\]
Then, writing $\bd$ for the differential of $K$, since
	\[\bd( t_j t_{j+1}t_{j+2}) = x_j t_{j+1}t_{j+2} - x_{j+1}t_j t_{j+2} + x_{j+2} t_j t_{j+1}
	\]
we deduce that there is no $a\in k$ such that $\bd( at_j t_{j+1}t_{j+2}) = x_{j+1}t_jt_{j+2}$. Hence $x_{j+1}t_j t_{j+2}$ is not a boundary in $K$, and so $\eta_{j,j+1,j+2}$ is not zero in $H_{2,\mbp_{j,3}}$.
\end{proof}
\end{prop}
 
Setting
	\[T = k\ang{\zeta_{1,2},\ldots,\zeta_{n-1,n},\eta_{1,2,3},\ldots,\eta_{n-2,n-1,n}},
	\]
there is a natural homomorphism $\varphi:T\to H$ of $\bbz^{n+1}$-graded algebras which sends the $\zeta_{i,i+1}$'s and $\eta_{j,j+1,j+2}$'s to themselves. By Propositions \ref{prop:BoocherEtAl}(1) and \ref{prop:bases}, it is surjective. Setting $I = \ker{\varphi}$, we thus have an isomorphism $T/ I\cong H'$ of graded algebras, and to prove $H$ is strand-Koszul it will suffice to show $I$ has a quadratic (with respect to strand degree) Gr\"obner basis.

Let $G$ be the subset of $T$ consisting of the following ten types of elements.  First, certain commutators and skew-commutators:
\begin{enumerate}
\item $\eta_{j,j+1,j+2}\eta_{i,i+1,i+2} - \eta_{i,i+1,i+2}\eta_{j,j+1,j+2}$ for all $1 \leq i<j \leq n-2$ (if $n\geq 4$),
\item $\eta_{j,j+1,j+2} \zeta_{i,i+1} - \zeta_{i,i+1} \eta_{j,j+1,j+2}$ for all $1 \leq i \leq j \leq n-2$,
\item $\zeta_{i,i+1} \eta_{j,j+1,j+2} - \eta_{j,j+1,j+2} \zeta_{i,i+1}$ for all $1\leq j < i \leq n-1$,
\item $\zeta_{j,j+1} \zeta_{i,i+1} + \zeta_{i,i+1} \zeta_{j,j+1}$ for all $1\leq i < j \leq n-1$.
\end{enumerate}
Then we have the ``string'' elements (so-called because the subscripts form strings of consecutive numbers):
\begin{enumerate}\setcounter{enumi}{4}
\item $\eta_{i,i+1,i+2}\zeta_{i+3,i+4} - \zeta_{i,i+1}\eta_{i+2,i+3,i+4} $ for all $1 \leq i \leq n-4$ (if $n\geq 5$),
\item $\zeta_{i,i+1}\zeta_{i+2,i+3}$ for all $1\leq i \leq n-3$ (if $n\geq 4$).
\end{enumerate}
Finally, the ``overlap'' elements:
\begin{enumerate}\setcounter{enumi}{6}
\item $\zeta_{i,i+1} \zeta_{j,j+1} $ for each $i=j-1,j$ with $1\leq i \leq j \leq n-1$,
\item $\zeta_{i,i+1} \eta_{j,j+1,j+2}$ for each $i=j-1,j,j+1,j+2$ with $1\leq i \leq n-1$ and $1\leq j \leq n-2$,
\item $\eta_{j,j+1,j+2}\zeta_{i,i+1}$ for each $i=j-1,j,j+1,j+2$ with $1\leq i \leq n-1$ and $1\leq j \leq n-2$, and
\item $\eta_{i,i+1,i+2}\eta_{j,j+1,j+2}$ for each $i=j-2,j-1,j$ with $1 \leq i \leq j \leq n-2$.
\end{enumerate}
Some elementary (if slightly tedious) computations show $(G)\subseteq I$.

We then consider the degree-lexicographic ordering on the monomials of $T$ with respect to the ordering 
	\[\zeta_{1,2} < \eta_{1,2,3} < \zeta_{2,3} < \eta_{2,3,4} <  \cdots < \zeta_{n-2,n-1} < \eta_{n-2,n-1,n} < \zeta_{n-1,n}.
	\]
The proof of the next lemma then easily follows.

\begin{lem}\label{lem:reduced}
Suppose $1\leq i,j \leq n-1$ and $1\leq \ell,m\leq n-2$. With respect to $G$,
\begin{enumerate}\label{lem:reduced}
\item the monomial $\zeta_{i,i+1}\zeta_{j,j+1}$ is reduced if and only if $i+1 \leq j-2 $;
\item the monomial $\zeta_{i,i+1} \eta_{\ell,\ell+1,\ell+2}$ is reduced if and only if $i+1 \leq \ell-1$;
\item the monomial $\eta_{\ell,\ell+1,\ell+2}\zeta_{i,i+1}$ is reduced if and only if $\ell+2 \leq i-2$;
\item the monomial $\eta_{\ell,\ell+1,\ell+2}\eta_{m,m+1,m+2}$ is reduced if and only if $\ell+2 \leq m-1$.
\end{enumerate}
\end{lem}

Now, suppose $1\leq i \leq n$, $1\leq r \leq n-i+1$, and $r\not \equiv 1 \modu{3}$. Define
	\[\mu_{i,r} = \begin{cases}
	\zeta_{i,i+1} & : r = 2, \\
	\zeta_{i,i+1} \eta_{i+2,i+3,i+4} \eta_{i+5,i+6,i+7} \cdots \eta_{i+r-3,i+r-2,i+r-1} & : r \equiv 2 \modu{3}, \ r>2, \\
	\eta_{i,i+1,i+2}\eta_{i+3,i+4,i+5} \cdots \eta_{i+r-3,i+r-2,i+r-1} & : r \equiv 0 \modu{3}.
	\end{cases}
	\]
By Lemma \ref{lem:reduced} each $\mu_{i,r}$ is reduced, and the following proposition shows that \textit{all} reduced monomials in $T$ are products of these monomials.

\begin{prop}\label{prop:charReduced}
Let $\mu\in T$ be a reduced monomial (not equal to $1$) of multigraded internal degree $\mbu\in \bbz^n$.
\begin{enumerate}
\item The multidegree $\mbu$ is squarefree.
\item If the complete decomposition of $\mbu$ is
	\[\mbu = \mbp_{i_1,r_1} + \cdots + \mbp_{i_m,r_m} \quad ( i_1 < \cdots < i_m),
	\]
then $r_p \not\equiv 1 \modu{3}$ for each $p=1,\ldots,m$ and $\mu = \mu_{i_1,r_1} \cdots \mu_{i_m,r_m}$.
\end{enumerate}
In particular, there is exactly one reduced monomial in $T_{\ast,\mbu}$.
\begin{proof}
We shall prove both statements simultaneously by inducing on the internal degree $|\mbu|$.

Note that each $\zeta_{i,i+1}$ has internal degree $2$ while each $\eta_{j,j+1,j+2}$ has internal degree $3$. Hence we must have $|\mbu|\geq 2$, and if equality holds, then necessarily $\mu = \zeta_{i,i+1} = \mu_{i,2}$ for some $i \in \{1,2,\ldots,n-1\}$. Similarly, if $|\mbu|=3$, then necessarily $\mu = \eta_{i,i+1,i+2} = \mu_{i,3}$ for some $i\in \{1,2,\ldots,n-2\}$.

Now assume $|\mbu|>3$ and that statements (1) and (2) both hold for reduced monomials of internal degree $<|\mbu|$. There are two cases to consider: either
	\[(\tn{a}) \ \mu = \nu \zeta_{i,i+1} \quad \tn{or} \quad (\tn{b}) \ \mu = \nu \eta_{i,i+1,i+2}
	\]
where $\nu$ is a submonomial of $\mu$ (not equal to $1$). In both cases since $\mu$ is reduced so is $\nu$, and our inductive hypotheses then imply that the multidegree $\mbv$ of $\nu$ is squarefree and that if $\mbv$ has complete decomposition
	\[\mbv = \mbp_{j_1,s_1} + \cdots + \mbp_{j_\ell,s_\ell} \quad (j_1 < \cdots < j_\ell),
	\]
then $s_q \not \equiv 1 \modu{3}$ for each $q=1,\ldots,\ell$ and $\nu = \mu_{j_1,s_1} \cdots \mu_{j_\ell,s_\ell}$. In case (a), since $\mu$ is reduced Lemma \ref{lem:reduced} implies $i\geq j_\ell + s_\ell+1$, which in turn implies that $\mbu$ (which is equal to $\mbv + \mbe_i + \mbe_{i+1}$) is squarefree with complete decomposition
	\[\mbu = \mbp_{j_1,s_1} + \cdots + \mbp_{j_\ell,s_\ell} + \mbp_{i,2}.
	\]
Since also $\mu = \mu_{j_1,s_1} \cdots \mu_{j_\ell,s_\ell} \mu_{i,2}$, we see that both statements (1) and (2) hold for $\mbu$ and $\mu$.

In case (b), Lemma \ref{lem:reduced} implies that $i\geq j_\ell+s_\ell$. If $i = j_\ell+s_\ell$, then $\mbu$ is squarefree with complete decomposition
	\[\mbu = \mbp_{j_1,s_1} + \cdots + \mbp_{j_\ell,s_\ell+3}
	\]
and $\mu = \mu_{j_1,s_1} \cdots \mu_{j_\ell,s_\ell+3}$. If $i > j_\ell+s_\ell$, then $\mbu$ is again squarefree with complete decomposition
	\[\mbu = \mbp_{j_1,s_1} + \cdots + \mbp_{j_\ell,s_\ell} + \mbp_{i,3}
	\]
and $\mu = \mu = \mu_{j_1,s_1} \cdots \mu_{j_\ell,s_\ell}\mu_{i,3}$. In either case, statements (1) and (2) both hold for $\mbu$ and $\mu$.
\end{proof}
\end{prop}

We now have enough to prove the main result.

\begin{thm}
Let $n$ be an integer $\geq 3$ and set
	\[R = \frac{k[X_1,X_2,\ldots,X_n]}{(X_1X_2,X_2X_3,\ldots,X_{n-1}X_n)}.
	\]
The Koszul homology algebra $H$ of $R$ is strand-Koszul.
\begin{proof}
We will prove that the natural surjection $\varphi:T \surj H$ induces an isomorphism $T/(G) \cong H$, having already shown $(G) \subseteq I$ where $I = \ker{\varphi}$.

For a fixed multidegree $\mbu\in \bbz^n$, the map $\varphi$ induces a surjection
	\[(T/(G))_{\ast,\mbu} \surj H_{\ast,\mbu}
	\]
of $k$-spaces. To show that the surjection is an isomorphism we will prove that
	\[\dim_k{(T/(G))_{\ast,\mbu}}= \dim_k{H_{\ast,\mbu}}. \tag{$\ast$}
	\]
Let $\mbu$ have complete decomposition $\sum_{p=1}^m \mbp_{i_p,r_p}$.

We recall that the images of the reduced monomials linearly span the quotient algebra $T/(G)$. Thus if $(T/(G))_{\ast,\mbu}\neq 0$, then by Proposition \ref{prop:charReduced} we have that $\dim_k{(T/(G))_{\ast,\mbu}}=1$, the multidegree $\mbu$ is squarefree, and $r_p \not \equiv 1 \modu{3}$ for each $p=1,\ldots,m$. By Proposition \ref{prop:BoocherEtAl}, we then have $\dim_k{H_{\ast,\mbu}}=1$. This establishes $(\ast)$.

Now, the isomorphism $T/(G)\cong H$ and Propositions \ref{prop:BoocherEtAl} and \ref{prop:charReduced} show that the images of the reduced monomials in $T/(G)$ form a basis, and hence the generators (1)-(10) of $(G)$ form a Gr\"obner basis. Each of these generators has strand degree $2$, and hence $H$ is strand-Koszul.
\end{proof}
\end{thm}

\begin{remark}
In addition to algebras cut out by edge ideals of paths, the paper \cite{BoocherEtAl2015} also contains results on algebras whose defining ideals are edge ideals of cycles. If $R$ denotes one of these latter types of algebras, then its Koszul homology algebra $H$ has a more complicated structure than what we've just seen above; indeed, Theorem 3.15 of \cite{BoocherEtAl2015} states that $H$ is no longer generated in its linear strand if $n \equiv 1 \modu{3}$ and $n>4$ (where $n$ is the number of vertices in the cycle). However, if $n\not\equiv 1 \modu{3}$, then just as we saw above, $H$ is generated by the subspaces $H_{1,2}$ and $H_{2,3}$; so, could it be that $H$ is strand-Koszul?

Not always. Indeed, if $n=9$, then running the following \texttt{Macaulay2} \cite{Macaulay2} code shows that $\beta^H_{2,6,9}(k)=1$, so that $H$ is not even strand quadratic:

\medskip

\begin{center}
\begin{minipage}{0.8\textwidth}
\texttt{needsPackage "DGAlgebras"}

\texttt{Q = QQ[$\mathtt{x_1}$..$\mathtt{x_9}$]}

\texttt{I = ideal($\mathtt{x_1}$*$\mathtt{x_2}$,$\mathtt{x_2}$*$\mathtt{x_3}$,$\mathtt{x_3}$*$\mathtt{x_4}$,$\mathtt{x_4}$*$\mathtt{x_5}$,$\mathtt{x_5}$*$\mathtt{x_6}$,$\mathtt{x_6}$*$\mathtt{x_7}$,$\mathtt{x_7}$*$\mathtt{x_8}$,$\mathtt{x_8}$*$\mathtt{x_9}$,$\mathtt{x_9}$*$\mathtt{x_1}$)}

\texttt{R = Q/I}

\texttt{H = HH koszulComplexDGA R}

\texttt{k = coker vars H}

\texttt{F = res(k,LengthLimit=>2)}

\texttt{peek betti F}
\end{minipage}
\end{center}
\end{remark}

%\vspace{-0.1in}
%%%%%%%%%%%%%%%%%%%%%%%%%%%%%%%%%%%%%%%%%%%%%%%
\bibliographystyle{amsplain}
\bibliography{/Users/johnmyers/GoogleDrive/tex/bib}

\end{document}